\documentclass[a4paper]{article}
\usepackage[latin1]{inputenc}
\usepackage{amsmath,amssymb}
\usepackage{amsfonts}
\newcounter{tnr}
\newcounter{crr}
\newenvironment{thm}%
{\medskip %
\bf   Theorem \arabic {tnr}. \it %
\stepcounter{tnr}%%
}%

\def\ff#1{{\cal F}_x^{#1}}
\def\mc#1{{\cal M}(#1)}
\def\smd{\sum_{d\leqslant xq_0q_1\atop d/(d,q_0q_2)\perp q_0q_1}}
\def\smn{\sum_{n\leqslant x\atop {n\perp q_0q_1 \atop d/(d,q_0q_2)|n}}}

\begin{document}
\title{On the densities of rational multiples}%Example 1:% Title 

 \author{Vilius Stak\. enas \\
 Faculty of Mathematics and Informatics\\
 Vilnius University\\
 email: vilius.stakenas@mif.vu.lt}

 \date{}
 \maketitle

 \begin{abstract}
For two subsets of natural numbers $ A,B\subset \mathbb{N} $ define the set of rational numbers $ \mc{A,B}$ with the elements 
represented by $ m/n, $ where $ m,n $ are coprime, $ m $ 
is divisible by some $ a\in A $  and $ n $ by some $ b\in B, $
respectively. Let $ I $ be some interval of positive real numbers and  $ {\cal F}_x^I $ denotes the set of rational 
numbers $ m/n \in I, $ such that $ m,n $ are coprime and 
$ n\leqslant x. $ The analogue to the Erd\"os-Davenport theorem about  multiples is proved: under some constraints on $ I $ the limits $ \sum\{\frac{1}{mn}:\frac{m}{n}\in {\cal F}_x^I\cap \mc{A,B} \}/ \sum\{\frac{1}{mn}:\frac{m}{n}\in {\cal F}_x^I\} $
exist for all subsets $ A,B\subset \mathbb{N}. $ 
 \end{abstract}

%
%================= Text entry area ================================%
\section{INTRODUCTION}\label{s:1}

For a subset $ A $ of natural numbers $ \mathbb{N} $ and $ x>1 $  denote 
\[ 
\nu_x^0(A)=\frac{1}{x} \sum_{n\in  A\cap [1;x]} 1 , \quad  \nu_x^1(A)=\frac{1}{\log x} \sum_{n\in  A\cap [1;x]} \frac{1}{n}.
 \]
The lower and upper limits as $ x\to \infty  $ will be denoted 
by $ \underline{\nu}^r(A), \overline{\nu}^r(A)\ (r=0,1); $   
the value of the limit if it  exists  by $ \nu^r(A) ,$ respectively.  

It follows from the chain of inequalities 
\[ 
\underline{\nu}^0(A)\leqslant \underline{\nu}^1(A)\leqslant 
\overline{\nu}^1(A)\leqslant \overline{\nu}^0(A) 
 \]
that the existence of $ \nu^0(A) $ implies the existence of 
$ \nu^1(A). $ If $ \nu^0(A) $ exists, we say that $ A $ possesses asymptotic density, and if  $ \nu^1(A) $ exists, $  A $ possesses 
logarithmic density. Even the  subsets $ A $ of  apparently simple structure may not possess asymptotic density. 

Let $ A\subset \mathbb{N}.$ The set of natural numbers divisible by some $ a\in A $ will be denoted by $\mc{A}, $ i.e. $ \mc{A} $
is the set of multiples of $ A. $

A.S. Besicovitch gave an example of $ A $ such that $ \mc{A} $ 
does not possess asymptotic density, see \cite{bes}. In 1937  H. Davenport and P. Erd\"os proved that every set of multiples have logarithmic density. Their original proof in \cite{derd} is based on Tauberian theorems, see also \cite{hall}, Theorem 02. The direct and elementary proof of this  theorem was provided  by the authors in \cite{dae}, it can be found also in the monograph  of H. Halberstam and K.F.  Roth, \cite{hroth}.  We   formulate the Erd\"os-Davenport theorem in the form, which results from 
the arguments in \cite{hroth}.

\begin{thm} Let $ A\subset \mathbb{N} $  and $ A_N=A\cap [1;N]$
for $ N\in \mathbb{N}. $
 Then $ \nu^1(\mc{A_N}), \nu^1(\mc{A}) $ exist, and
\begin{equation}
\nu^1(\mc{A}) =\lim_{N\to \infty }\nu^1(\mc{A_N}). \label{tankn}
 \end{equation} 
\end{thm}

The main aim of this paper is to investigate the density questions related to the sets of multiples of rational numbers. 

Let  $  \mathbb{Q}^+ $ be  the set of positive rational numbers. For the natural numbers $ m,n $ we denote as usually by $(m,n)$
their greatest common divisor.   If $ (m,n)=1, $ i.e. the numbers are coprime, we write $ m\perp n $ (suggestion of R.L. Graham,  D.E. Knuth and O. Potashnik, see \cite{grah}, p.115). For the rational numbers $ r\in  \mathbb{Q}^+ $ we shall always use the unique representation $ r=m/n, m,n\in \mathbb{N}, m\perp n. $

For two subsets $ A,B\subset \mathbb{N} $ and $ q\in \mathbb{N} $ we define the set of multiples in $ \mathbb{Q}^+ $ by
\[ 
\mc{A,B|q}=\Big\{\frac{m}{n}:m\in \mc {A}, n\in \mc {B}, mn \perp q \Big\}.
 \] 
If $ q=1 $ we write $ \mc{A,B}$  instead of  $\mc{A,B|1}. $

Let $ I_x=(\lambda_1(x),\lambda_2(x)) $ be some system of intervals, $ I_x\subset (0; +\infty ), x\geqslant 1. $ We shall write in the following briefly $ I=(\lambda_1, \lambda_2) $ and introduce the sets of rational numbers 
\[ \ff{I}=\Big\{\frac{m}{n}: \frac{m}{n}\in \mathbb{Q}^+, n\leqslant x \Big\}\cap I. \]  
Let $ R\subset \mathbb{Q}^+ $ and $ r_1,r_2 \in \{0,1\}.$ Then if $\ff{I}\not =\emptyset  ,$ we denote
\[ 
S_{x,I}^{r_1r_2}(R)=\sum_{m/n \in \ff{I}\cap R} m^{-r_1}n^{-r_2}, \quad \nu_x^{r_1r_2}(R)=\frac{S_{x,I}^{r_1r_2}(R)}{S_{x,I}^{r_1r_2}(\ff{I})}.
 \] 
 
 If the limit of $\nu_x^{r_1r_2}(R)  $ exists for $ R\subset \mathbb{Q}^+ $ as $ x\to \infty , $ it will be denoted by $\nu^{r_1r_2}(R),$ and the lower and upper limits  by $ \underline{\nu}^{r_1r_2}(R), \overline{\nu}^{r_1r_2}(R), $ respectively.

We investigate the limit behaviour of $\nu_x^{r_1r_2}(\mc{A,B|q})  $   as $ x\to \infty  $ under some conditions 
imposed on $ \lambda_i. $ In the case of unit interval $ I=(0,1) $ related problems were  considered in authors paper \cite{plng}.

\section{OVERVIEW OF RESULTS}\label{s:2}

If interval $ I=(\lambda_1,\lambda_2) $ does not depend on $ x, $ the inequalities of type (\ref{tankn}) can be proved.

\begin{thm} Let the interval $ I=(\lambda_1, \lambda_2) $ be fixed. Then for an arbitrary $ A\subset \mathbb{N} $
\begin{eqnarray*} 
&&\underline{\nu}^{00}(A)\leqslant  \underline{\nu}^{01}(A)\leqslant \overline{\nu}^{01}(A)\leqslant \overline{\nu}^{00}(A),\\
&&\underline{\nu}^{10}(A)\leqslant  \underline{\nu}^{11}(A)\leqslant \overline{\nu}^{11}(A)\leqslant \overline{\nu}^{10}(A).
\end{eqnarray*}
\end{thm}

If $ A,B $ are finite subsets of $ \mathbb{N} $  the following statement holds.

\begin{thm} Let $ \lambda_1<\lambda_2 $ satisfy the following conditions:

if $ \lambda_1=0, $ then $ \lambda_2>x^{-c} $ for some $ 0<c<1; $

if $ \lambda_1>0,$ then $ \lambda_1\log(\lambda_2/\lambda_1)\log x \to \infty  $ as $ x\to \infty . $

Then for finite sets $ A,B\subset \mathbb{N} $ and $ q\in \mathbb{N} $ all densities 
$ \nu^{r_1r_2}\big(\mc{A,B|q} \big) $ exist and are equal.

\end{thm}

\medskip

Note that if $ \lambda_1>0 $ and $ (\lambda_2-\lambda_1)/\lambda_1 $ remains bounded, 
the constraints on $ \lambda_i$ are equivalent to requirement
$(\lambda_2-\lambda_1)\cdot \log x \to \infty   $ as $ x\to \infty . $

It is possible to prove under appropriate conditions on $ \lambda_i $ this statement for the sets satisfying
\[ \sum_{d\in A\cup B}\frac{1}{d}<\infty , \]
 but we shall not pursue this question.

The inequality for densities in the following theorem 
should be compared to Heilbronn-Rohrbach inequality proved in 
\cite{heilb}, \cite{rohr}; see also \cite{hall}.

\begin{thm} Let the sets  $ A,B\subset \mathbb{N} $
be finite  and satisfy the following conditions: $ a\perp b $ for all $ a\in A, b\in B; $ if $ a_1,a_2\in A, b_1,b_2\in B, $
then $ a_1\perp a_2/(a_1,a_2), b_1\perp b_2/(b_1,b_2). $ Let  $ \nu\big(\mc{A,B|q}\big) $ denote the common value 
of densities from Theorem 3. Then  the following inequality holds:
\[ 
1-\nu \big(\mc{A,B|q}\big) \geqslant \prod_{p|q}\Big(1-\frac{2}{p+1}\Big)\cdot   \prod_{c\in A\cup B}\Big( 
1-\frac{1}{c}\prod_{p|c}\Big(1-\frac{1}{p+1} \Big)\Big).
 \]
\end{thm} 

The sets satisfying conditions of Theorem 4 can be constructed as follows. Let $ r_1,r_2,\ldots  $ be an arbitrary sequence of coprime integers. If $ a_j=\prod_{k\in I_j } r_k, $ where
$ I_j $ is some finite subset of naturals then, obviously, 
$ a_i\perp a_j/(a_j,a_i) $ for all pairs $ i,j. $

The main result of the paper is an analogue or Erd\"os-Davenport 
theorem for the sets of rational multiples.

\begin{thm} Let for  the intervals $I=(\lambda_1,\lambda_2)  $
the following conditions be satisfied: 

if $ \lambda_1=0 $ then $ \lambda_2 x \to \infty  $ and
 $ \log x/\log(\lambda_2x)<c_1 $ as $ x\to \infty  $ with some $ c_1>0; $

if $ \lambda_1>0 $ then with some positive constants $ c_2,c_3 $
\[ 
\frac{c_2}{\log(\lambda_2+2)}<\lambda_1<\lambda_2<x^{c_3}, \quad
\frac{1}{\log (\lambda_2+2)}\cdot \log\Big(\frac{\lambda_2}{\lambda_1}\Big)\cdot \log x \to \infty , \quad x\to \infty.
 \]
  Then  for arbitrary $ A,B\subset \mathbb{N} $  and $ q\in \mathbb{N} $ the limit
\[ \nu^{11}\big(\mc{A,B|q} \big)=\lim_{x\to \infty }\nu_x^{11}\big(\mc{A,B|q} \big) \] 
exists.
\end{thm} 
 
Let $ \lambda_1>c (c>0) $ and $ \lambda_2 $ be bounded. Then the conditions of Theorem 5 for $ \lambda_i $ can be reduced to requirement
\[(\lambda_2-\lambda_1)\cdot \log x \to \infty   \quad \text{as} \quad  x\to \infty . \]
 
 \section{PROOFS}\label{s:3}

 Let $ q_0,q_1,q_2 $ be some coprime natural numbers and 
 \begin{equation}
\mathbb{Q}_{q_0,q_1,q_2}=\Big\{\frac{m}{n}\in \mathbb{Q}^+, mn\perp q_0, mq_1\perp nq_2 \Big\}. \label{qq}
  \end{equation}
We investigate the asymptotical behaviour of the  sums $ S_{x,I}^{r_1r_2}(\mathbb{Q}_{q_0,q_1,q_2}) $ as $ x\to \infty . $ Methods beeing used are elementary, the remainder terms in the asymptotics depend on $ q_i. $  

{\bf Lemma. \it  Let for the coprime integers $ q_0,q_1,q_2 $
\[
\Pi(q_0,q_1,q_2)=\prod_{p|q_0}\Big(1-\frac{2}{p+1}\Big)
\prod_{p|q_1q_2}\Big(1-\frac{1}{p+1}\Big).
  \]
 Then the following asymptotics hold 
 
\begin{eqnarray*}
\frac{S^{00}_{x,I}(\mathbb{Q}_{q_0,q_1,q_2})}
{\Pi(q_0,q_1,q_2)}&=&\frac{3}{\pi^2}(\lambda_2-\lambda_1)x^2\Big\{1+O\Big(\frac{\log x}{x}+\frac{\log x}{(\lambda_2-\lambda_1)x} \Big)\Big\},\\
\frac{S^{01}_{x,I}(\mathbb{Q}_{q_0,q_1,q_2})}{\Pi(q_0,q_1,q_2)}&=&\frac{6}{\pi^2}(\lambda_2-\lambda_1)x\Big\{1+O\Big(\frac{\log x}{x}+\frac{\log^2 x}{(\lambda_2-\lambda_1)x} \Big)\Big\}.
\end{eqnarray*} 
 If $ \lambda_1>0 $ then
 \begin{eqnarray*}
 \frac{S^{10}_{x,I}(\mathbb{Q}_{q_0,q_1,q_2})}{\Pi(q_0,q_1,q_2)}&=&\frac{6}{\pi^2}\log\Big(\frac{\lambda_2}{\lambda_1}\Big) x\Big\{1+O\Big(\frac{\log x}{x}+\frac{\log^2 x}{\lambda_1 \log\big(\frac{\lambda_2}{\lambda_1}\big)x} \Big)\Big\},\\
\frac{S^{11}_{x,I}(\mathbb{Q}_{q_0,q_1,q_2})}{\Pi(q_0,q_1,q_2)}&=&\frac{6}{\pi^2}\log\Big(\frac{\lambda_2}{\lambda_1}\Big)\log x\Big\{1+O\Big(\frac{1}{\log x}+\frac{1}{\lambda_1\log\big(\frac{\lambda_2}{\lambda_1}\big)\log x} \Big) \Big\}. 
 \end{eqnarray*}
In the case $ \lambda_1=0 $ we have
\begin{eqnarray*}
\frac{S^{10}_{x,I}(\mathbb{Q}_{q_0,q_1,q_2})}{\Pi(q_0,q_1,q_2)}&=&\frac{6}{\pi^2}x\log(\lambda_2x)
\Big\{1+O\Big( \frac{1}{\log(\lambda_2x)}+\frac{\log x}{x} \Big) \Big \},\\
\frac{S^{11}_{x,I}(\mathbb{Q}_{q_0,q_1,q_2})}{\Pi(q_0,q_1,q_2)}&=&\begin{cases}
\frac{3}{\pi^2}\log^2(\lambda_2x)\Big\{1+O\Big(\frac{\log x}{\log^2(\lambda_2x)} \Big) \Big\}, &\text{if $ \frac{1}{x}<\lambda_2\leqslant  1, $}\\
\frac{3}{\pi^2}\log x\cdot \log(\lambda_2^2x)\Big\{1+O\Big(\frac{\log (\lambda_2x)}{\log x\cdot \log(\lambda_2^2x)} \Big) \Big\}, &\text{if $ \lambda_2>1. $} \label{qvv}\\
\end{cases}
\end{eqnarray*} 
} 

The functions in O-signs of the Lemma are diferrent. It is easily seen, that if $ \lambda_1=0, $ then the condition
$ x^{-c}<\lambda_2 $ with some $ 0<c<1 $ is sufficient for all functions in O-signs related to the case $ \lambda_1=0 $ 
to be vanishing. 

Consider now the case $ \lambda_1>0. $ The function
\[ 
f(u)=u-c\log\big(1+\frac{u}{c}\big), \quad u\geqslant 0, \quad c>0,
 \]
 is not decreasing, hence 
 \[ 
\lambda_1\log\Big(\frac{\lambda_2}{\lambda_1} \Big)=
\lambda_1\log\Big(1+\frac{\lambda_2-\lambda_1}{\lambda_1}\Big)\leqslant \lambda_2-\lambda_1. 
  \]
It follows from this, that under condition
\[ 
\lambda_1\log\Big(\frac{\lambda_2}{\lambda_1} \Big) \frac{x}{\log^2 x} \to \infty , \quad  x\to \infty ,
 \]
all functions in O-signs of $ S^{r_1r_2}_{x,I}(\mathbb{Q}_{q_0,q_1,q_2}),$  with  $r_1+r_2<2, $ 
are vanishing. We include $ S_{x,I}^{11} $ if we use  the stronger  requirement 
\[ 
\lambda_1\log\Big(\frac{\lambda_2}{\lambda_1} \Big) \log x \to \infty , \quad  x\to \infty .
 \] 
If $ q_0=q_1=q_1=1, $ then
\[ 
S^{r_1r_2}_{x,I}(\mathbb{Q}_{1,1,1})=S^{r_1r_2}_{x,I}(\mathbb{Q}^+)=
\sum\Big\{m^{-r_1}n^{-r_2}:\frac{m}{n}\in \ff{I}\Big\}.
 \] 
The following Corollary  follows easily from the Lemma. 

{\bf Corollary. \it 
Let $ \lambda_i $ fulfill the following conditions
\begin{eqnarray*}
&&\text{if } \lambda_1=0\ \text{then } x^{-c}<\lambda_2 \text{ with some } 0<c<1;\\
&&\text{if } \lambda_1>0 \ \text{then } \lambda_1\log\Big(\frac{\lambda_2}{\lambda_1} \Big) \log x \to \infty , \quad  x\to \infty  .
\end{eqnarray*}

Then for all $ r_1,r_2 $ and fixed coprime numbers $ q_0,q_1,q_2 $
 \[ 
\frac{S^{r_1r_2}_{x,I}(\mathbb{Q}_{q_0,q_1,q_2})}{S^{r_1r_2}_{x,I}(\mathbb{Q}^+)} \to \Pi(q_0,q_1,q_2) \quad as \quad x\to \infty . 
  \]
}

{\bf Proof.}  We abbreviate the notation as $ S^{r_1r_2}=S^{r_1r_2}_{x,I}(\mathbb{Q}_{q_0,q_1,q_2}) $
and start with the expression
\[ 
S^{r_1r_2}=\sum_{n\leqslant x \atop n\perp q_0q_1}n^{-r_2}\sum_{\lambda_1n<m<\lambda_2n \atop m\perp nq_0q_2}m^{-r_1}.
 \] 
With the M\"obius function  $ \mu(n) $ we proceed as follows
\begin{eqnarray}
S^{r_1r_2}&=&\sum_{n\leqslant x \atop n\perp q_0q_1}n^{-r_2}\sum_{\lambda_1n<m<\lambda_2n }m^{-r_1}\sum_{d|(m,nq_0q_2)}\mu(d)=
\sum_{d\leqslant xq_0q_1}\mu(d)\sum_{n\leqslant x\atop {n\perp q_0q_1 \atop d|nq_0q_2}}n^{-r_2}\sum_{\lambda_1n<m<\lambda_2n \atop d|m}m^{-r_1} \nonumber \\ 
&=&\smd \mu(d)d^{-r_1}\sum_{n\leqslant x\atop {n\perp q_0q_1 \atop d/(d,q_0q_2)|n}}n^{-r_2}\sum_{\lambda_1\frac{n}{d}<m<\lambda_2\frac{n}{d}}m^{-r_1}. \label{pr0}
\end{eqnarray}

For the last sum over $ m $  we shall use the folowing equalities 
\[ 
\sum_{\lambda_1\frac{n}{d}<m<\lambda_2\frac{n}{d}}m^{-r_1}=
\begin{cases}(\lambda_2-\lambda_1)\frac{n}{d}+\theta_{n,d},&\text{if $ r_1=0 $},\\
\log\big(\frac{\lambda_2}{\lambda_1}\big)+\theta_{n,d}\frac{d}{\lambda_1n},& \text{if $ \lambda_1>0, r_1=1 $},\\
\log\big(\frac{\lambda_2n}{d}\big)+\theta_{n,d},& \text{if $ \lambda_1=0, r_1=1, $ and $ \frac{\lambda_2n}{d}>1, $}
\end{cases}
 \]
 where $ \theta_{n,d} $ are bounded by some absolute constant. 
 
 Consider  the case $ r_1=r_2=0 $ first. Then 
 \begin{equation}
S^{00}=(\lambda_2-\lambda_1) \smd\frac{\mu(d)}{d}
\sum_{n\leq x\atop {n\perp q_0q_1 \atop d/(d,q_0q_2)|n}}n+
O\Big(\sum_{d\leqslant xq_0q_1}\mu^2(d)\sum_{n\leqslant x\atop {n\perp q_0q_1 \atop d/(d,q_0q_2)|n}}1 \Big). \label{pr1}
\end{equation}

Let $ S^{00}_1 $  stands for the main term in (\ref{pr1}). Using the divisibility property $ d/(d,q_0q_2)|n $ and the asymptotics 
\[ 
\sum_{n\leqslant u \atop n\perp q}n=\frac{1}{2}u^2\prod_{p|q}\Big(1-\frac{1}{p}\Big) +O(u),
 \]  
 we rewrite the main  term of $ S^{00} $ as
\begin{equation}
S^{00}_1=\frac{1}{2}(\lambda_2-\lambda_1)x^2\prod_{p|q_0q_1}\Big(1-\frac{1}{p}\Big)\smd(d,q_0q_2)\frac{\mu(d)}{d^2}+O\big((\lambda_2-\lambda_1)x\log x\big).\label{pr2}
 \end{equation} 
Note that $ d/(d,q_0q_2)\perp q_0q_1 $ is equivalent to $ d\perp q_1, $ hence  
\begin{eqnarray*}
\smd (d,q_0q_2)\frac{\mu(d)}{d^2}&=&\sum_{d\perp q_1}(d,q_0q_2)\frac{\mu(d)}{d^2}+O\Big(q_0q_2\sum_{d>xq_0q_1}\frac{1}{d^2}\Big)\\
&=&\prod_{p}\Big(1-\frac{1}{p^2}
 \Big)\prod_{p|q_1}\Big(1-\frac{1}{p^2} \Big)^{-1}
\prod_{p|q_0q_2}\Big(1+\frac{1}{p} \Big)^{-1}+O\Big(\frac{1}{x}\Big).
 \end{eqnarray*}
Setting this in  (\ref{pr2}) one gets
\[ 
S^{00}_1=\frac{3}{\pi^2}(\lambda_2-\lambda_1)x^2\Pi(q_0,q_1,q_2)+O\big((\lambda_2-\lambda_1)x\log x\big).
 \]
For the remainder term  in (\ref{pr1}) we use the bound  
\[ 
\sum_{d\leqslant xq_0q_1}\mu^2(d)\sum_{n\leqslant x\atop {n\perp q_0q_1 \atop d/(d,q_0q_2)|n}}1\leqslant 
\sum_{d\leqslant xq_0q_1}\mu^2(d)\sum_{n\leqslant x(d,q_0q_2)/d}1=O\big(x\log x\Big).
 \]
Hence putting all together we obtain
\[ 
S^{00}=\frac{3}{\pi^2}(\lambda_2-\lambda_1)x^2\Pi(q_0,q_1,q_2)\Big\{1+O\Big(\frac{\log x}{x}+\frac{\log x}{(\lambda_2-\lambda_1)x} \Big)\Big\}.
 \] 
Consider now the case $ r_1=0, r_2=1. $  Then instead of (\ref{pr1}) we have
 \begin{equation}
S^{01}=(\lambda_2-\lambda_1) \smd\frac{\mu(d)}{d}
\sum_{n\leqslant x\atop {n\perp q_0q_1 \atop d/(d,q_0q_2)|n}}1+
O\Big(\sum_{d\leqslant xq_0q_1}\mu^2(d)\sum_{n\leqslant x\atop {n\perp q_0q_1 \atop d/(d,q_0q_2)|n}}\frac{1}{n} \Big). \label{pr3}
\end{equation}
Let  $ S^{01}_1 $ denote the main term in (\ref{pr3}). Using 
\[ 
\sum_{n\leqslant u \atop n\perp q}1=u\prod_{p|q}\Big(1-\frac{1}{p}\Big) +O(1)
 \]  
we obtain
\begin{eqnarray*}
S^{01}_1&=&(\lambda_2-\lambda_1)x\prod_{p|q_0q_1}\Big(1-\frac{1}{p}\Big)\smd(d,q_0q_2)\frac{\mu(d)}{d^2}+O\big((\lambda_2-\lambda_1)\log x\big)\\ 
&=&\frac{6}{\pi^2}(\lambda_2-\lambda_1)x\Pi(q_0,q_1,q_2)+O\big((\lambda_2-\lambda_1)\log x\big).
\end{eqnarray*}
For the remainder term in (\ref{pr3}) use the obvious bound
\[ 
\sum_{d\leqslant xq_0q_1}\mu^2(d)\sum_{n\leqslant x\atop {n\perp q_0q_1 \atop d/(d,q_0q_2)|n}}\frac{1}{n}\leqslant
\sum_{d\leqslant xq_0q_1}(d,q_0q_2)\frac{\mu^2(d)}{d}
\sum_{n\leqslant x}\frac{1}{n}=O\big(\log^2x \big).
 \]
 Hence the asymptotics 
\[ 
S^{01}=\frac{6}{\pi^2}(\lambda_2-\lambda_1)x\Pi(q_0,q_1,q_2)\Big\{1+O\Big(\frac{\log x}{x}+\frac{\log^2 x}{(\lambda_2-\lambda_1)x} \Big)\Big\}.
 \] 
 is established.
 
Suppose now that $ r_1=1,r_2=0. $ From  (\ref{pr0}) one gets
\[ 
S^{10}=\smd \frac{\mu(d)}{d}\smn \sum_{\lambda_1 \frac{n}{d} <m<\lambda_2 \frac{n}{d}} \frac{1}{m}.
 \] 
Let first $ \lambda_1>0. $ Then
 \[ 
S^{10}=\log \Big(\frac{\lambda_2}{\lambda_1}\Big) \smd \frac{\mu(d)}{d}\smn 1+O\Big(\frac{1}{\lambda_1}\sum_{d\leqslant xq_0q_1}\mu^2(d)\smn \frac{1}{n}\Big).
  \]
Expression for $ S^{10} $ differs from that one in (\ref{pr3})
in term involving $ \lambda_i $ only. Hence, in the same way as above we get
 \[ 
S^{10}=\frac{6}{\pi^2}\log\Big(\frac{\lambda_2}{\lambda_1}\Big) x\Pi(q_0,q_1,q_2)\Big\{1+O\Big(\frac{\log x}{x}+\frac{\log^2 x}{\lambda_1 \log\big(\frac{\lambda_2}{\lambda_1}\big)x} \Big)\Big\}.
 \] 
 Let now $ \lambda_1=0. $ Then 
 \[ 
S^{10}=\smd \frac{\mu(d)}{d}\smn \log \Big(\frac{\lambda_2n}{d}\Big) +O\Big(\sum_{d\leqslant xq_0q_1}\frac{\mu^2(d)}{d}\smn 1\Big).
  \]
  The remainder term does not exceed 
\[ 
x\sum_{d}(d,q_0q_2)\frac{\mu^2(d)}{d^2}=O(x).  
\]
Using the divisibility condition $ d/(d,q_0q_2)|n $ we proceed as follows
\begin{eqnarray*}
S^{10}&=&\smd \frac{\mu(d)}{d}\sum_{n\leqslant x(d,q_0q_2)/d \atop 
n\perp q_0q_1} \log \Big(\frac{\lambda_2n}{(d,q_0q_2)}\Big)+O(x)
=
\smd \frac{\mu(d)}{d}\sum_{n\leqslant x(d,q_0q_2)/d \atop 
n\perp q_0q_1} \log (\lambda_2n)\\ &-&
\smd \frac{\mu(d)}{d}\log(d,q_0q_2)\sum_{n\leqslant x(d,q_0q_2)/d \atop 
n\perp q_0q_1}1+O(x).
 \end{eqnarray*}
The second minus term is $ O(x), $ hence
\begin{equation}
S^{10}=\smd \frac{\mu(d)}{d}\sum_{n\leqslant x(d,q_0q_2)/d \atop 
n\perp q_0q_1} \log (\lambda_2n)+O(x). \label{pr4}
 \end{equation}
 Using 
 \[ 
\sum_{n\leqslant u \atop n\perp q}1=u\prod_{p|q}\Big(1-\frac{1}{p}\Big)+O(1) 
  \]
 and integrating by parts one derives for $ c>0 $ easily 
 \[ 
 \sum_{n\leqslant u \atop n\perp q}\log(cn)=u\log(cu)\prod_{p|q}\Big(1-\frac{1}{p}\Big)+O(u+|\log(cu)|), \quad \text{as} \quad u\to \infty .
  \] 
   
  Using this in (\ref{pr4}) we get
\begin{eqnarray*} 
S^{10}&=&\prod_{p|q_0q_2}\Big(1-\frac{1}{p} \Big)x
\smd\frac{\mu(d)}{d^2}(d,q_0q_2)\log \Big(\lambda_2x\frac{(d,q_0q_2}{d} \Big)\\
&+&
O\Big(x+x\smd\frac{\mu^2(d)}{d^2}(d,q_0q_2)+
\smd\frac{\mu^2(d)}{d}\Big|\log \Big(\lambda_2x\frac{(d,q_0q_2}{d} \Big)\Big| \Big).
 \end{eqnarray*}
 It is easily seen that the remainder term can be reduced 
 to $ O(x+\log x\cdot \log(\lambda_2x)+\log^2x)=O(x+\log x \log(\lambda_2x^2)). $ Using 
 additivity property for the logarithm in the first sum 
 we split the main term of $ S^{10} $ into two parts and the second will be $ O(x). $ Hence 
 \[ 
S^{10}= \prod_{p|q_0q_2}\Big(1-\frac{1}{p} \Big)x \log(\lambda_2x)
\smd\frac{\mu(d)}{d^2}(d,q_0q_2)+O(x+\log x\cdot \log(\lambda_2x^2)).
  \]
 The remaining sum was calculated above, then simplifying the remainder terms one gets
\[ 
S^{10}=\frac{6}{\pi^2}\Pi(q_0,q_1,q_2)x\log(\lambda_2x)
\Big\{1+O\Big( \frac{1}{\log(\lambda_2x)}+\frac{\log x}{x} \Big) \Big \}.
 \]  
 With $ r_1=r_2=1 $ we have 
 \[ 
S^{11}=\smd\frac{\mu(d)}{d}\smn\frac{1}{n}\sum_{\lambda_1\frac{n}{d}<m<\lambda_2\frac{n}{d}}\frac{1}{m}. 
  \]
If $ \lambda_1>0 $ this reduces to
\[ 
S^{11}=\log\Big(\frac{\lambda_2}{\lambda_1}\Big)\smd\frac{\mu(d)}{d}\smn\frac{1}{n}+O\Big(\frac{1}{\lambda_1}\smd\frac{\mu^2(d)}{d}\cdot d \smn\frac{1}{n^2} \Big).
 \]  
 The sum over $ n $ in the remainder term is $ O(d^{-2}), $ hence
 \begin{eqnarray*} 
S^{11}&=&\log\Big(\frac{\lambda_2}{\lambda_1}\Big)\smd\frac{\mu(d)}{d}\smn\frac{1}{n}+O(\lambda_1^{-1})\\
&=&\log\Big(\frac{\lambda_2}{\lambda_1}\Big)\smd\frac{\mu(d)}{d^2}(d,q_0q_1)\sum_{n\leqslant x(d,q_0q_2)/d\atop n\perp q_0q_1}\frac{1}{n}+O(\lambda_1^{-1}). 
  \end{eqnarray*}
Using the asymptotics
\[ 
\sum_{n\leqslant u \atop n\perp q}\frac{1}{n}=\prod_{p|q}\Big(
1-\frac{1}{p}\Big)\log u+O(1),
 \]  
 we derive
 \begin{eqnarray*}
S^{11}&=&\log\Big(\frac{\lambda_2}{\lambda_1}\Big)\prod_{p|q_0q_1}\Big(1-\frac{1}{p}\Big)\smd\frac{\mu(d)}{d^2}(d,q_0q_2)\log\Big(x\cdot \frac{(d,q_0q_2)}{d} \Big)+O\Big(\frac{1}{\lambda_1}+\log\Big(\frac{\lambda_2}{\lambda_1}\Big)\Big)\\
&=&
\log\Big(\frac{\lambda_2}{\lambda_1}\Big)\prod_{p|q_0q_1}\Big(1-\frac{1}{p}\Big)\log x\smd\frac{\mu(d)}{d^2}(d,q_0q_2)+O\Big(\frac{1}{\lambda_1}+\log\Big(\frac{\lambda_2}{\lambda_1}\Big)\Big).
\end{eqnarray*}
Simplifying the sum over $ d $ as above we arrive finally to 
\[ 
S^{11}=\frac{6}{\pi^2}\log\Big(\frac{\lambda_2}{\lambda_1}\Big)\Pi(q_0,q_1,q_2)\log x\Big\{1+O\Big(\frac{1}{\log x}+\frac{1}{\lambda_1\log\big(\frac{\lambda_2}{\lambda_1}\big)\log x} \Big) \Big\}.
 \]  
 
Consider now the case $ \lambda_1=0: $
\begin{eqnarray*}
S^{11}&=&\smd\frac{\mu(d)}{d}\sum_{d/\lambda_2<n\leqslant x
\atop{n\perp q_0q_1 \atop d/(d,q_0q_2)|n}}\frac{1}{n}\sum_{m<\lambda_2\frac{n}{d}}\frac{1}{m}\\
&=&\smd\frac{\mu(d)}{d}
\sum_{d/\lambda_2<n\leqslant x
\atop{n\perp q_0q_1 \atop d/(d,q_0q_2)|n}}
 \frac{1}{n}\log\Big(\frac{\lambda_2n}{d} \Big)+
O\Big(\smd \frac{\mu^2(d)}{d}\smn \frac{1}{n} \Big).
\end{eqnarray*}
Using the divisibility condition $d/(d,q_0q_2)|n  $ we reduce the  term in O-sign to  $ O(\log x) $ and simplify the expression as follows
\begin{eqnarray*} 
S^{11}&=&\smd\frac{\mu(d)}{d^2}(d,q_0q_2)
\sum_{(d,q_0q_2)/\lambda_2<n\leqslant x (d,q_0q_2)/d \atop
n\perp q_0q_1}\frac{1}{n}\log\Big(\frac{\lambda_2n}{(d,q_0q_2)} \Big)+O(\log x)\\
&=&\smd\frac{\mu(d)}{d^2}(d,q_0q_2)
\sum_{(d,q_0q_2)/\lambda_2<n\leqslant x (d,q_0q_2)/d \atop
n\perp q_0q_1}\frac{\log(\lambda_2n)}{n}+O(\log x).
\end{eqnarray*}
Extending the sum over $ n $ to the range $ 1/\lambda_2<n\leqslant x $ we introduce the error term $ O(\log x +\log(\lambda_2x)). $ Hence 
\[ 
S^{11}=\smd\frac{\mu(d)}{d^2}(d,q_0q_2)
\sum_{1/\lambda_2<n\leqslant x \atop
n\perp q_0q_1}\frac{\log(\lambda_2n)}{n}+
O(\log x +\log(\lambda_2x)).
 \]
The main term is expressed as the product of two sums,
the first one equals to
\[ \frac{6}{\pi ^2}\Pi(q_0,q_1,q_2)\prod_{p|q_0q_1}\Big(1-\frac{1}{p}\Big)^{-1}+O(x^{-1}). \] 
The second sum of the main term can  be calculated by partial integration, the final result would be 
\[ 
\sum_{1/\lambda_2<n\leqslant x\atop n\perp q_0q_1}\frac{\log\big(\lambda_2n \big)}{n} =
\begin{cases}
\frac{1}{2}\prod_{p|q_0q_1}\Big(1-\frac{1}{p}\Big)(\log x+\log\lambda_2)^2+O(1),&\text{as $ \lambda_2<1, $}\\
\frac{1}{2}\prod_{p|q_0q_1}\Big(1-\frac{1}{p}\Big)(\log^2 x+2\log\lambda_2\log x)+O(\log(\lambda_2+1)),&\text{as $ \lambda_2\geqslant 1 $}.
\end{cases}
 \] 
If we write $ (\log x+\log\lambda_2)^2=\log^2(\lambda_2x) $
and $ (\log^2 x+2\log\lambda_2\log x)=\log x \cdot \log(\lambda_2^2x), $ then after manipulating with the remainder terms we arrive to the following expressions
\[ 
S^{11}=\begin{cases}
\frac{3}{\pi^2}\log^2(\lambda_2x)\Big\{1+O\Big(\frac{\log x}{\log^2(\lambda_2x)} \Big) \Big\}, &\text{if $ \frac{1}{x}<\lambda_2\leqslant 1, $}\\
\frac{3}{\pi^2}\log x\cdot \log(\lambda_2^2x)\Big\{1+O\Big(\frac{\log (\lambda_2x)}{\log x\cdot \log(\lambda_2^2x)} \Big) \Big\}, &\text{if $ \lambda_2>1. $}\\
\end{cases}
 \]
Note that the remainder term for $ \lambda_2>x^{-c}, $ where
$ 0<c<1, $ is  $ O(\log^{-1} x). $ 
The Lemma is proved.

{\bf Proof of Theorem 2.} Let us start with the first chain of inequalities. Because of the interval $  I $ is fixed 
\[ 
S^{00}_{x,I}(\mathbb{Q}^+)\sim \frac{3}{\pi^2}|I|x^2, \quad 
S^{01}_{x,I}(\mathbb{Q}^+)\sim \frac{6}{\pi^2}|I|x, \quad x\to \infty,
 \]
 where $ |I|=\lambda_2-\lambda_1. $ For an arbitrary subset
 $ A\subset \mathbb{Q}^+ $ we have 
 \begin{equation}
 S^{01}_{x,I}(A)=\int_{1-}^x\frac{1}{t}{\rm d}S^{00}_{t,I}(A)=
 \frac{S^{00}_{x,I}(A)}{x}+\int_{1-}^x\frac{S^{00}_{t,I}(A)}{t^2}{\rm d}t .\label{intnl}
 \end{equation}
For an arbitrary fixed $ \epsilon>0 $ we shall have
\[ 
S^{00}_{t,I}(A)\leqslant (\overline{\nu}^{00}(A)+\epsilon )\frac{3}{\pi^2}|I|t^2
 \]
 as $ t\geqslant t_0. $ From this observation and (\ref{intnl})
 we derive 
 \[ 
 S^{00}_{x,I}(A)\leqslant (\overline{\nu}^{00}(A)+\epsilon )\frac{6}{\pi^2}|I|x+C, 
  \]
  with some $ C>0. $ Then, consequently, $ \overline{\nu}^{01}(A)
  \leqslant \overline{\nu}^{00}(A). $ The inequality for lower limits follows   from the inequality for complement set $ \overline{\nu}^{01}(A^c)
  \leqslant \overline{\nu}^{00}(A^c). $

The second chain of inequalities can be derived in an analogous manner  from the equality
\[ 
S^{11}_{x,I}(A)=\int_{1-}^x\frac{1}{t}{\rm d}S^{10}_{t,I}(A).
 \]
Theorem 2 is proved.  

{\bf Proof of Theorem 3.} 
Consider now the sets of multiples $ \mc{A,B|q}. $ If $ A=\{a \},
B=\{b\} $ we shall write $ \mc{A,B|q}=\mc{a,b|q}. $
For natural numbers $ a,b $  with $ (a,b)>1 $ or $ (ab,q)>1 $ we have  $ \mc{a,b|q}=\emptyset . $ Let  $ a\perp b   $  and 
$ ab\perp q. $  Then using the notation (\ref{qq})
\[ 
S^{r_1r_2}_{x,I}(\mc{a,b|q})=a^{-r_1}b^{-r_2}S^{r_1r_2}_{b^{-1}x,ba^{-1}I}(\mathbb{Q}_{q,a,b}).
 \]
After examining  the asymptotics of Lemma  we conclude that  under  conditions  of Theorem 3 for $ \lambda_i  $
\[ 
S^{r_1r_2}_{b^{-1}x,ba^{-1}I}(\mathbb{Q}_{q,a,b}) \sim a^{r_1-1}
b^{r_2-1}S^{r_1r_2}_{x,I}(\mathbb{Q}_{q,a,b}), \quad x\to \infty . \]
From the Corollary  we obtain
\begin{equation} 
\nu^{r_1r_2}_x(\mc{a,b|q}) \to \frac{1}{ab}
\prod_{p|q}\Big(1-\frac{2}{p+1} \Big)
\prod_{p|ab}\Big(1-\frac{1}{p+1} \Big), \quad x\to \infty .\label{tank}
 \end{equation} 
 
Let now $ A,B $ be two finite sets. By the sieve arguments we have
\begin{equation} 
\nu_x^{r_1r_2}(\mc{A,B|q})=\sum_{C\subset A\times B \atop C\not =\emptyset }(-1)^{|C|+1}\nu_x^{r_1r_2}\Big(\bigcap_{\langle a, b\rangle \in C}\mc{a,b|q}\Big). \label{ret}
 \end{equation}
For $ C=\{\langle a_1,b_1\rangle , \langle a_2,b_2\rangle , \ldots ,\langle a_{|C|},b_{|C|}\rangle  \}\subset A\times B $ 
let us introduce the notations
\[ 
[C]_A=[a_1,a_2,\ldots , a_{|C|}], \quad [C]_B=[b_1,b_2,\ldots , b_{|C|}],
 \]
here $ [...] $ stands for the least common multiples of numbers in the brackets.  Then clearly
\[ 
\bigcap_{\langle a, b\rangle \in C}\mc{a,b|q}=\mc{[C]_A,[C]_B|q}.
 \]
Due to (\ref{tank}) all the summands in (\ref{ret}) tend to their limits as $ x\to \infty . $ Hence the statement of Theorem 3 follows.

{\bf Proof of Theorem 4.} The inequality follows by induction over the number of elements $ |A|+|B|. $ If $ A=\{a\},   B=\{b\}, $ then  either $ \nu\big(\mc{a,b|q} \big)=0 $ or 
\begin{equation}
\nu\big(\mc{a,b|q} \big)=\prod_{p|q}\Big(1-\frac{2}{p+1} \Big)
\frac{1}{ab}\prod_{p|ab}\Big(1-\frac{1}{p+1} \Big).
 \end{equation} \label{denst}
In the first case the inequality is trivial, and in the second one we have 
\[ 
1-\nu\big(\mc{a,b|q} \big)=1-\prod_{p|q}\Big(1-\frac{2}{p+1} \Big)
\prod_{p|ab}\Big(1-\frac{1}{p+1} \Big)\geqslant \prod_{p|q}\Big(1-\frac{2}{p+1} \Big) \Big(1-\prod_{p|ab}\Big(1-\frac{1}{p+1} \Big)\Big).
 \] 
 Let the inequality holds for some finite sets $ A,B $ and we 
 add a new number $ a^* $ to $ A. $ We shall show that the inequality will be satisfied for $ \mc{A^*,B|q} $ with 
 $ A^*=A\cup \{a^*\}, $ too.
 Let us introduce the following notations: $ [a^*,A]=\{[a^*,a]:a\in A \}, ]a^*,A[=\{a/(a,a^*):a\in A\}, $
 where $ [a^*,a] $ denotes the least common multiple of numbers in brackets; if $ C $ is some finite set of numbers, then $ [C] $ stands for the least common multiple of all elements of $ C. $ We start with 
 \[ 
\mc{A^*,B|q}=\mc{A,B|q}\cup  \big(\mc{A^*,B|q}\backslash 
\mc{A,B|q}\big).
  \]
Denote briefly $ \mc{A,B|q}^*=\mc{A^*,B|q}\backslash 
\mc{A,B|q}. $ Then 
\begin{eqnarray*}
\nu\big(\mc{A^*,B|q} \big)&=&\nu\big(\mc{A,B|q} \big)+
\nu\big(\mc{A,B|q}^* \big),\\
\nu\big(\mc{A,B|q}^* \big)&=&\nu\big(\mc{a^*,B|q}\big)-\nu\big(\mc{[a^*,A],B|q} \big).
\end{eqnarray*}
Using the sieve arguments  and the properties of $ A $ one derives
\begin{eqnarray*}
\nu\big(\mc{a^*,B|q}\big)&=&\sum_{C\subset B \atop C\not = \emptyset}(-1)^{1+|C|}  \nu\big(\mc{a^*,[C]|q} \big)=
\frac{1}{a^*}\prod_{p|a^*}\Big(1-\frac{1}{p+1} \Big)\nu\big( \mc{1,B|q}\big),\\
\nu\big(\mc{[a^*,A],B|q}\big)&=&\sum_{C\subset [a^*,A]\times B \atop C\not = \emptyset}(-1)^{1+|C|}  \nu\big(\mc{[C]_{[a^*,A]},[C]_B|q} \big)\\ &=&
\frac{1}{a^*}\prod_{p|a^*}\Big(1-\frac{1}{p+1} \Big)\nu\big( \mc{]a^*,A[,B|q}\big).
\end{eqnarray*} 
It follows now from this that 
\[ 
\nu\big(\mc{A^*,B|q} \big)=\nu\big(\mc{A,B|q}\big)+
\frac{1}{a^*}\prod_{p|a^*}\Big(1-\frac{1}{p+1} \Big)\big(\nu\big(\mc{1,B|q} \big)-\nu\big(\mc{]a^*,A[,B|q }\big) \big).
 \] 
Because of $\nu\big(\mc{1,B|q} \big)\leqslant 1  $ and 
$\nu\big(\mc{]a^*,A[,B|q }\big) \geqslant \nu\big(\mc{A,B|q }\big)  $ 
we obtain
\[ 
1-\nu\big(\mc{A^*,B|q} \big) \geqslant 1-\nu\big(\mc{A,B|q} \big) - \frac{1}{a^*}\prod_{p|a^*}\Big(1-\frac{1}{p+1} \Big)\big(1-\nu\big(\mc{A,B|q }\big) \big),
 \] 
 and the inequality for the sets $ A^*, B $ follows. If instead of $ A $ we  add a new element to $ B, $ the arguments proving the inequality would be essentially the same. The Theorem is proved. 
 
{\bf Proof of Theorem 5.} 
Recall that for $ A\subset \mathbb{N} $ we denote by $ \mc{A} $ the set of multiples of elements $ a\in A. $  If  $ N> 1 $  let 
$ A_N=A\cap [1;N]. $

We  start with the equality
\begin{equation}
\nu_x^{r_1r_2}(\mc{A,B|q}) = \nu_x^{r_1r_2}(\mc{A_N,B_N|q})+
\nu_x^{r_1r_2}(\mc{A,B|q}\backslash \mc{A_N,B_N|q}). \label{tnk11}
 \end{equation} 
 
It suffices to show that for any $ \epsilon>0 $  the upper limit of the second term in 
(\ref{tnk11}) is less than $ \epsilon $ as 
$ x\to \infty , $ supposed that $ N $ is large enough.
Define two subsets of rational numbers 
\[
{\cal M}_N^1=\Big\{\frac{m}{n}: m\in {\cal M}(A)\backslash {\cal M}(A_N) \Big\}, \quad 
{\cal M}_N^2=\Big\{\frac{m}{n}: n\in {\cal M}(B)\backslash {\cal M}(B_N) \Big\} .
\]
Then 
\[ 
\mc{A,B|q}\backslash \mc{A_N,B_N|q}\subset 
{\cal M}_N^1\cup {\cal M}_N^2.
 \]
We are going to prove that for fixed $ \delta >0 $ and  $ N $ sufficiently large we shall have $ \overline{\nu }^{11}\big( 
{\cal M}_N^i \big)\leqslant \delta $ for $ i=1,2. $ Denote for the sake of brevity  
$ \mc{A}_N= \mc{A}\backslash \mc{A_N},
\mc{B}_N= \mc{B}\backslash \mc{B_N}.
$
Then 
\[ 
S^{11}_{x,I}\big({\cal M}_N^2 \big)\leqslant 
\sum_{n\leqslant x\atop n\in\mc{B}_N }\frac{1}{n}\sum_{\lambda_1n<m<\lambda_2 n}\frac{1}{m}.
 \]
Let $ \lambda_1=0 $ first.  With some constant $ c>0 $ we have 
\[ 
S^{11}_{x,I}\big({\cal M}_N^2 \big)\leqslant (\log(\lambda_2x)+c)\sum_{n\leqslant x\atop n\in\mc{B}_N }\frac{1}{n}.
 \]
The Erd\"os-Davenport statement as formulated in the Theorem 1 implies that there exists some vanishing sequence $ \delta_N $
such that  
$ \overline{\nu}^1\big(\mc{B}_N \big)<\delta_N . $ It follows then  that for $ x  $ a sufficiently large we shall have
\[ 
S^{11}_{x,I}\big({\cal M}_N^2 \big)\leqslant \delta_N \log x (\log(\lambda_2x)+c).
 \]
Compare now the functions on the right-side of this inequality to that ones in the asymptotics of $ S^{11}_{x,I}(\mathbb{Q}^+) $
(see Lemma):
\[ 
S^{11}_{x,I}\big({\cal M}_N^2 \big)\leqslant 
\begin{cases} \log^2(\lambda_2x)\big\{\delta_N  \frac{\log x}{\log(\lambda_2x)}+\frac{\delta_N c\log x}{\log^2(\lambda_2x)}\big\},& \\
\log x\cdot \log(\lambda_2^2x)\big\{\delta_N \frac{\log(\lambda_2x)}{\log(\lambda_2^2x)}+\frac{\delta_N c}{\log(\lambda_2^2x)}\big\}.&  
\end{cases}
 \]
Having in mind the conditions on $ \lambda_i $ we conclude that $   \overline{\nu}^{11}\big({\cal M}_N^2 \big)\leqslant \delta  $ 
for $ N $ large enough. 

We shall show now that $\overline{\nu}^{11}\big({\cal M}_N^1 \big) \leqslant \delta $
as well. If $ m/n<\lambda_2 $ and $ n\leqslant x, $ then $ m\leqslant \lambda_2 x $ and $ n>m/\lambda_2. $
We start with 
\[ 
S^{11}_{x,I}\big({\cal M}_N^1 \big)\leqslant \sum_{m\leqslant \lambda_2 \atop m\in\mc{A}_N  }
\frac{1}{m}\sum_{m/\lambda_2 <n\leqslant x}\frac{1}{n} +\sum_{\lambda_2<m\leqslant \lambda_2x \atop m\in\mc{A}_N  }
\frac{1}{m}\sum_{m/\lambda_2 <n\leqslant x}\frac{1}{n}.
 \]
 Consider the first summand. Using the Erd\"os-Davenport theorem as above we obtain that for $ x $ large enough  
 \[ 
\sum_{m\leqslant \lambda_2 \atop m\in\mc{A}_N  }
\frac{1}{m}\sum_{m/\lambda_2 <n\leqslant x}\frac{1}{n}\leqslant (\log x+c)\sum_{m\leqslant \lambda_2 \atop m\in\mc{A}_N  }
\frac{1}{m} \leqslant \delta_N (\log x+c)\log \lambda_2. 
  \]
where $ \delta_N \to 0 $ as $ N\to \infty . $
Using similar arguments for the second sum we get
\begin{eqnarray*}
S^{11}_{x,I}\big({\cal M}_N^1 \big)&\leqslant & \delta_N (\log x+c)\log \lambda_2+ 
\sum_{\lambda_2<m\leqslant \lambda_2x \atop m\in\mc{A}_N  }
\frac{1}{m}\Big\{\log\Big(\frac{\lambda_2x}{m}\Big)+c\frac{\lambda_2}{m}  \Big\}\\ 
&\leqslant & \delta_N (\log x+c)\log \lambda_2+
\log x\sum_{m\leqslant \lambda_2x \atop m\in\mc{A}_N  }\frac{1}{m}+c\lambda_2\sum_{m>\lambda_2}\frac{1}{m^2}\leqslant 
\delta_N \log x \log (\lambda_2 x) +c_1, \quad c_1>0,
 \end{eqnarray*}
and $ \overline{\nu}^{11}\big({\cal M}_N^1 \big) \leqslant \delta . $ This completes the proof in the case $ \lambda_1=0. $

Let now $ \lambda_1>0. $ Then using the Erd\"os-Davenport theorem again we have 
\begin{eqnarray*} 
S^{11}_{x,I}\big({\cal M}_N^2 \big)&\leqslant & 
\sum_{n\leqslant x\atop n\in\mc{B}_N }\frac{1}{n}\sum_{\lambda_1n<m<\lambda_2 n}\frac{1}{m}\leqslant 
\sum_{n\leqslant x\atop n\in\mc{B}_N }\frac{1}{n}
\Big\{\log\Big(\frac{\lambda_2}{\lambda_1} \Big)+\frac{1}{\lambda_1 n} \Big\}\\
&\leqslant & \log\Big(\frac{\lambda_2}{\lambda_1} \Big) \log x\Big(\delta_N +\frac{1}{\lambda_1\log(\lambda_2/\lambda_1)\log x}\Big).
\end{eqnarray*}
Note that under conditions on $ \lambda_1, \lambda_2 $
\[
\lambda_1\log(\lambda_2/\lambda_1)\log x =\lambda_1\log(\lambda_2+2)\cdot \Big(\frac{1}{\log (\lambda_2+2)} \log \Big( \frac{\lambda_2}{\lambda_1}\Big) \log x \Big) \to \infty ,
 \]
as $ x\to \infty . $ Consequently $ \overline{\nu}^{11}\big({\cal M}_N^2 \big) \leqslant \delta . $

For $ S^{11}_{x,I}\big({\cal M}_N^1 \big) $ we proceed as follows:
\begin{equation} 
S^{11}_{x,I}\big({\cal M}_N^1 \big)\leqslant \sum_{m\leqslant \lambda_2 \atop m\in \mc{A}_N}\frac{1}{m} \sum_{n<m/\lambda_1}\frac{1}{n}+
\sum_{\lambda_2<m\leqslant \lambda_2x \atop m\in \mc{A}_N} \sum_{m/\lambda_2<n<m/\lambda_1}\frac{1}{n}. \label{s11}
 \end{equation} 

If $ \lambda_2 $ remains bounded, then the first sum in (\ref{s11}) is zero for $  N $ sufficiently large. Otherwise we have
\[ 
\sum_{m\leqslant \lambda_2 \atop m\in \mc{A}_N}\frac{1}{m} \sum_{n<m/\lambda_1}\frac{1}{n} \leqslant 
\sum_{m\leqslant \lambda_2 \atop m\in \mc{A}_N} \frac{1}{m}\Big\{ 
\log\Big(\frac{\lambda_2}{\lambda_1}\Big)+c\Big\} \leqslant 
\log\Big(\frac{\lambda_2}{\lambda_1}\Big)\log x\Big\{\delta_N 
\frac{\log \lambda_2}{\log x}+\delta_N \frac{c\log \lambda_2}{\log(\lambda_2/\lambda_1)\log x} \Big\}.
 \] 
For  the second sum  in (\ref{s11}) we obtain
\begin{eqnarray*}
\sum_{\lambda_2<m\leqslant \lambda_2x \atop m\in \mc{A}_N} \frac{1}{m} \sum_{m/\lambda_2<n<m/\lambda_1}\frac{1}{n} &\leqslant & \sum_{\lambda_2<m \leqslant  \lambda_2x \atop m\in \mc{A}_N} \frac{1}{m} 
\Big\{ 
\log\Big(\frac{\lambda_2}{\lambda_1}\Big)+c\frac{\lambda_2}{m} \Big\} 
\leqslant  
\delta_N \log\Big(\frac{\lambda_2}{\lambda_1}\Big) \log(\lambda_2x)+c \\ &\leqslant & 
\log\Big(\frac{\lambda_2}{\lambda_1}\Big)\log x\Big\{\delta_N 
\frac{\log (\lambda_2x)}{\log x}+\delta_N \frac{c}{\log(\lambda_2/\lambda_1)\log x} \Big\}.
 \end{eqnarray*}
It follows from both estimates that for given $ \delta >0 $ under conditions on $ \lambda_i $  we shall have 
$S^{11}_{x,I}\big({\cal M}_N^1 \big)\leqslant  \delta 
S^{11}_{x,I}\big(\mathbb{Q}^+\big), $
supposed $ x,N $ are large enough. Hence  
 $ \overline{\nu}^{11}\big({\cal M}_N^1 \big) \leqslant \delta , $ and the proof of theorem is completed.

\bibliographystyle{plain}

%spell_to this point %============== End of text entry ===============%

\bibliography{x}

\end{document}